\documentclass{elsart3-1}

\usepackage{amssymb}
\usepackage[T1]{fontenc}
\usepackage[latin1]{inputenc}
\usepackage{latexsym}
\usepackage{xypic}
\usepackage{eufrak}
\usepackage{euscript}
\usepackage{amsmath}
\usepackage{verbatim}
\usepackage{fancyhdr}

\usepackage[english,francais]{babel}

\newtheorem{theoreme}{Th\'eor\`eme}[section]
\newtheorem{lemm}[theoreme]{Lemme}

\newtheorem{defi}[theoreme]{D\'efinition\rm}

\def\og{\leavevmode\raise.3ex\hbox{$\scriptscriptstyle\langle\!\langle$~}}
\def\fg{\leavevmode\raise.3ex\hbox{~$\!\scriptscriptstyle\,\rangle\!\rangle$}}

\def\Cu{\EuScript{C}^1}

\def\Nd{\EuScript{N}^1}

\renewcommand{\d}{\displaystyle}

\newcommand{\K}{\mathbb{K}}
\newcommand{\N}{\mathbb{N}}

\newcommand{\R}{\mathbb{R}}

\newcommand{\C}{\mathbb{C}}

\def\cal{\mathcal}

\begin{document}
% Vous pouvez mettre dans la prochain ligne la rubrique choisie
% (si vous la connaissez) - meme deux, format : Rubrique1/Rubrique2
\centerline{Probability Theory/Dynamical Systems}
\begin{frontmatter}

% Titre, auteurs et adresses

% utiliser la commande \thanksref dans \title, \author ou \address
%     pour les notes en bas de page ;

% utiliser la commande \ead pour l'adresse e-mail de chaque auteur
%    (apr\"{E}s la commande \auteur) ;

% \author{Name\thanksref{label2}}
\selectlanguage{francais}
\title{\textbf{Lemme de coh\'{e}rence et th\'{e}or\`{e}me de Noether
stochastique}}

% utiliser les \`{E}tiquettes pour indiquer l'adresse de chaque auteur,
%     s'il y a plusieurs adresses


\author[cresson]{Jacky CRESSON},
\ead{cresson@math.univ-fcomte.fr}
\author[darses]{S\'{e}bastien DARSES}
\ead{darses@math.univ-fcomte.fr}

\address[cresson]{Laboratoire de Math\'{e}matiques, Universit\'{e} de Franche-Comt\'{e}}
\address[darses]{Laboratoire de Math\'{e}matiques, Universit\'{e} de Franche-Comt\'{e}}
% etc, etc

% Vous pouvez mettre a la prochaine ligne les dates
% (de reception et d'acceptation), et le nom du presentateur de votre Note

\medskip
\selectlanguage{francais}
%\begin{center}
%{\small Re\c{c}u le *****~; accept\'e apr\`es r\'evision le +++++\\
%Pr\'esent\'e par {\pounds}{\pounds}{\pounds}{\pounds}{\pounds}}
%\end{center}

\begin{abstract}
% resume en francais, et apres l'abstract en anglais, qui
%    commence avec le titre en gras.
\selectlanguage{francais}

La proc\'{e}dure de plongement stochastique, d\'{e}finie dans
\cite{note1}, permet d'associer \`{a} l'\'{e}quation d'Euler-Lagrange
classique (EL) une \'{e}quation d'Euler-Lagrange stochastique (ELS).
Cette derni\`{e}re est-elle sous-tendue par un principe de moindre
action g\'{e}n\'{e}ralis\'{e} ? Pour aborder cette question, nous d\'{e}veloppons
un calcul des variations stochastique, initi\'{e} par Yasue
\cite{ya1}. On donne un analogue stochastique naturel $F$ de la
fonctionnelle lagrangienne d'action. On introduit une notion de
stationnarit\'{e} pour laquelle les solutions de (ELS) sont les points
stationnaires de $F$. La notion de stationnarit\'{e} ainsi d\'{e}finie
rend coh\'{e}rent le calcul des variations stochastique vis-\`{a}-vis de
la proc\'{e}dure de plongement stochastique. Enfin, nous d\'{e}montrons un
th\'{e}or\`{e}me de Noether stochastique qui sugg\`{e}re d'introduire une
nouvelle notion, celle d'int\'{e}grale premi\`{e}re stochastique.

{\it Pour citer cet article~: A. Nom1, A. Nom2, C. R. Acad. Sci.
Paris, Ser. I 340 (2005).} \vskip 0.5\baselineskip

\selectlanguage{english} \noindent{\bf Abstract} \vskip
0.5\baselineskip \noindent {\bf Coherence lemma and stochastic
Noether theorem. } The stochastic embedding procedure defined in
\cite{note1} associates a stochastic Euler-Lagrange equation (SEL)
to the standard Euler-Lagrange equation (EL). Can we derive (SEL)
from a generalized least action principle? To address this
question, we develop a stochastic calculus of variation initiated
by Yasue \cite{ya1}. We give a stochastic analog $F$ of the
lagrangian action functional. We introduce a notion of
stationarity according to which the solutions of (SEL) are the
stationary points of $F$. This notion of stationarity brings
coherence to stochastic calculus of variation with respect to
stochastic embedding. Finally, we prove a stochastic Noether
theorem which introduces an original notion of stochastic first
integral.
 {\it To cite this article: A. Nom1, A. Nom2, C. R.
Acad. Sci. Paris, Ser. I 340 (2005).}
\end{abstract}
\end{frontmatter}

% Maintenant la version abr\`{E}g\`{E}e en anglais, si pr\`{E}sente

%\selectlanguage{english}
%\section*{Abridged English version}

% Texte de la version abr\`{E}g\`{E}e en anglais

\selectlanguage{francais}

\section{D\'efinition d'un calcul des variations stochastique}
\label{}

On note $I:=]a,b[$ o\`u $a<b$ et $J:=[a,b]$ l'adh\'erence de $I$
dans $\R$. Soit $\K$ un corps et $d\in\N^*$. On se donne un espace
probabilis\'e $(\Omega,\cal{A},P)$ sur lequel existent une famille
croissante de tribus $(\cal{P}_t)_{t\in J}$ et une famille
d\'ecroissante de tribus $(\cal{F}_t)_{t\in J}$.

\begin{defi}

On note $\Cu_{\K}(J)$ l'ensemble des processus $X$ d\'efinis sur
$J\times \Omega$, \`a valeurs dans $\K^d$ et tels que : $X$ soit
$(\cal{P}_t)$ et $(\cal{F}_t)$ adapt\'{e}, pour tout $t\in J$ $X_t\in
L^2(\Omega)$, l'application $t\to X_t$ de $J$ dans $L^2(\Omega)$
est continue, pour tout $t\in I$ les quantit\'es
$DX_t=\lim_{h\rightarrow 0^+}h^{-1} E[X_{t+h}-X_t\mid {\cal P}_t
]$, et $D_* X_t=\lim_{h\rightarrow 0^+} h^{-1} E[X_t-X_{t-h}\mid
{\cal F}_t ]$,
existent dans $L^2(\Omega)$, et enfin les applications $t\to DX_t$ et $t\to D_*X_t$ sont continues de $I$ dans $L^2(\Omega)$.\\
Le compl\'et\'e de $\Cu_{\K} (J)$ pour la norme $\parallel
X\parallel=\sup_{t\in I} (\parallel X_t\parallel_{L^2(\Omega)}
+\parallel DX_t\parallel_{L^2(\Omega)} +\parallel D_*
X(t)\parallel_{L^2(\Omega)} )$, est encore not\'e $\Cu_{\K}(J)$,
et simplement $\Cu(J)$ quand $\K=\R$.
\end{defi}

On note $\cal{D}$ la d\'{e}riv\'{e}e stochastique introduite dans
\cite{stoc} et d\'{e}finie par ${\cal D} =\d {D+D_* \over 2} +i {D-D_*
\over 2}$.

On dira qu'un lagrangien $L$ est admissible si la fonction
$L(x,v)$ est d\'{e}finie sur $\R^d\times\C^d$, $C^1$ en $x$ et
holomorphe en $v$, et est r\'{e}elle quand $v$ est r\'{e}el. $L$ est dit
naturel s'il s'\'{e}crit $L(x,v)=q(v)-U(x)$ o\`{u} $q$ est une forme
quadratique sur $\C^d$ et $U$ un potentiel de classe $\cal{C}^1$
sur $\R^d$. Soit alors la fonctionnelle associ\'{e}e \`{a} $L$ d\'{e}finie par
$F_J : \Xi\subset\Cu(J)\to \C$, $\d
F_J(X)=E\left[\int_JL(X_t,\cal{D}X_t)dt\right]$ avec $\d
\Xi=\left\{X\in\Cu(J),E\left[\int_J
|L(X_t,\cal{D}X_t)|dt\right]<\infty\right\}$. On d\'{e}finit
l'ensemble $\cal L$ des processus $L$-adapt\'{e}
$\cal{L}=\left\{X\in\Cu(J),\partial_x L(X_t,\cal{D}X_t)\in\Cu(J),
\partial_v L(X_t,\cal{D}X_t)\in\Cu(J)\right\}$.

Soit $\Gamma$ un sous espace de $\Cu(J)$. On appelle
$\Gamma-$variation d'un processus $X\in\Cu(J)$, un processus de la
forme $X+Z$ o\`{u} $Z\in\Gamma$. Soit $\Gamma_{\Xi}$ le sous-espace
vectoriel de $\Gamma$ d\'{e}fini par $\d
\Gamma_{\Xi}=\left\{Z\in\Gamma, \forall X\in\Xi,
Z+X\in\Xi\right\}$, et $\Nd(J)$ le sous-espace vectoriel
$\Nd(J)=\{X\in \Cu(J), DX=D_*X\}$. Posons alors la

\begin{defi}
Si $L$ est un Lagrangien admissible et $F_J$ la fonctionnelle
associ\'{e}e, $F_J$ est dite $\Gamma$-diff\'{e}rentiable en un processus
$X\in\Xi$  si pour tout $\d Z\in\Gamma_{\Xi}$,
$F_J(X+Z)-F_J(X)=dF_J (X,Z)+R_X(Z)$ , o\`{u} $dF_J(X,Z)$ est une
fonctionnelle lin\'{e}aire en $Z\in \Gamma_{\Xi}$ et
$R_X(Z)=o(\parallel Z\parallel )$. De plus $X$ est dit
$\Gamma-$stationnaire si pour tout $\d Z\in\Gamma_{\Xi}$,
$dF_J(X,Z)=0$.
\end{defi}

On consid\`{e}re le cas $\Gamma=\Nd(J)$. La $\Nd(J)$-diff\'{e}rentielle de
$F_J$ est donn\'{e}e par:

\begin{lemm} Soit $L$ un lagrangien admissible dont toutes les diff\'{e}rentielles secondes sont born\'{e}es.
Posons $g(Z,\partial_v L)(s)=E\left[Z_s\partial_v
L(X_s,\mathcal{D}X_s)\right]$. Alors $\Nd (I)_{\Xi}=\Nd(I)$ et la
fonctionnelle $F_J$ associ\'{e}e \`{a} $L$ est $\Nd (I)$-diff\'{e}rentiable en
tout processus $\d X\in\Xi\cap\cal{L}$, et pour tout $\d Z\in\Nd
(I)$, sa
diff\'{e}rentielle s'\'{e}crit :\\
$\d dF_J(X,Z) = E\left [\int_a^b\ (\partial_x L-\cal{D}\partial_v
L)(X_u,\mathcal{D}X_u) Z_u du \right ] +g(Z,\partial_v
L)(b)-g(Z,\partial_v L)(a)$.

\end{lemm}

\textbf{D\'{e}monstration.} \`A l'aide du d\'{e}veloppement de Taylor de
$L$, on a $\d L(X+Z,\mathcal{D}(X+Z))-L(X,\mathcal{D}X) =
\partial_x L(X,\mathcal{D}X)Z +\partial_v
L(X,\mathcal{D}X)\mathcal{D}Z+\int_0^1 (1-t)\left(\partial_x^2
L(T^t)Z^2+\partial_{xv}^2L(T^t)Z\cal{D}Z+\partial_v^2
L(T^t)(\cal{D}Z)^2\right)dt$ o\`{u} $T^t=(X+tZ,\cal{D}X+t\cal{D}Z)$.
D'apr\`{e}s l'in\'{e}galit\'{e} de Cauchy-Schwarz,  $\sup_J E[|Z\cal{D}Z|]$,
$\sup_J E[|Z|^2]$ et  $\sup_J E[|\cal{D}Z|^2]$ sont des
$O(\|Z\|^2)$ car $Z\in\Cu(J)$. De plus $\partial_{xv}^2L$,
$\partial_{x}^2L$ et $\partial_{v}^2L$ sont born\'{e}es et $\d
X\in\Xi\cap\cal{L}$. On en d\'{e}duit que $\Nd (I)_{\Xi}=\Nd(I)$ et
$\d F_J(X+Z)-F_J(X)=E\left[\int_a^b\left(\partial_x
L(X_s,\mathcal{D}X_s)Z_s +\partial_v
L(X_s,\mathcal{D}X_s)\mathcal{D}Z_s\right)ds\right]+o(\|Z\|)$. On
conclut en utilisant (\ref{formule2}) et le fait que $Z\in\Nd(J)$.
$\square$

Rappelons les lemmes suivants respectivement d\'{e}montr\'{e}s dans
\cite{stoc} p.42 et \cite{stoc} p.70, le premier g\'{e}n\'{e}ralisant la
"loi produit" donn\'{e}e par Nelson dans \cite{ne1} p.$80$, le
deuxi\`{e}me \'{e}tant l'analogue stochastique du lemme classique des
"fonctions plateaux" (\textit{cf} \cite{ar} p.$57$) :
\begin{lemm}\label{formule2}
Soit $X,Y\in\EuScript{C}^1_{\C}(I)$. Alors $\d
E[\mathcal{D}X_t\cdot Y_t+X_t\cdot \overline{\mathcal{D}}Y_t] =
\frac{d}{dt}E[X_t\cdot Y_t]$.
\end{lemm}

\begin{lemm}
Soit $Y\in\Cu(J)$. Si pour tout $Z\in\Nd(J)$ $\d \int_J E[Y_u
\cal{D}Z_u]du=0$, alors $Y$ est constant.
\end{lemm}

Des trois lemmes on d\'{e}duit le

\begin{theoreme}[$\Nd$-Principe de moindre action stochastique]
Une condition n\'{e}cessaire et suffisante pour qu'un processus
$X\in\Xi\cap\cal{L}$ soit un processus $\Nd(J)-$stationnaire pour
$F_J$ est qu'il v\'{e}rifie l'\'{e}quation d'Euler-Langrange stochastique
$(ELS)$ : $(\partial_x L-\cal{D}\partial_v
L)(X_u,\mathcal{D}X_u)=0$ sur $J$.

\end{theoreme}

Le calcul des variations stochastiques est d\'{e}velopp\'{e} ici
ind\'{e}pendemment de la proc\'{e}dure de plongement stochastique d\'{e}finie
dans \cite{note1}. Il devient coh\'{e}rent avec cette derni\`{e}re pour un
choix ad\'{e}quat de l'espace de variation, en ce sens que :

\begin{lemm}[Lemme de coh\'{e}rence]
Le diagramme suivant commute :
\begin{eqnarray}
\xymatrix{
  & L(x(t),x'(t)) \ar[d]_{\bigstar} \ar[r]^{\cal{S}} & L(X_t,\mathcal{D}X_t)
  \ar[d]^{ \maltese}       \\
  & (EL) \ar[r]_{\cal{S}}   & (ELS)
  }
\end{eqnarray}
o\`{u} $\bigstar$ symbolise l'utilisation du principe de moindre
action classique \`{a} partir de la fonctionnelle\\
$\d \int_J L(x(t),x'(t))dt$, $\maltese$ celle du $\Nd$-principe de
moindre action stochastique \`{a} partir de $\d
E\left[\int_JL(X_t,\cal{D}X_t)dt\right]$ et $\cal S$ la proc\'{e}dure
de plongement stochastique d\'{e}finie dans \cite{note1}.
\end{lemm}

Dans \cite{stoc} chap.$7$, on traite le cas $\Gamma=\Cu(J)$ et on
montre qu'une condition n\'{e}cessaire et suffisante pour qu'un
processus $X\in\Xi\cap\cal{L}$ soit un processus
$\Cu(J)-$stationaire pour $F_J$ est qu'il v\'{e}rifie l'\'{e}quation :
$(\partial_x L-\overline{\cal{D}}\partial_v
L)(X_u,\mathcal{D}X_u)=0$ sur $J$. Ce cas ne nous permet pas
d'obtenir un lemme de coh\'{e}rence mais un th\'{e}or\`{e}me de Noether
stochastique.

\section{Th\'{e}or\`{e}me de Noether stochastique}

Les sym\'{e}tries apparaissant dans certains syst\`{e}mes lagrangiens
induisent l'existence d'int\'{e}grales premi\`{e}res du mouvement
(\textit{cf} \cite{ar} p.$88$). Une question naturelle est alors
d'\'{e}tudier la persistance de ces objets fondamentaux sur le syst\`{e}me
stochastis\'{e} par la proc\'{e}dure de plongement. Introduisons
l'ensemble $\cal P$ des processus d\'{e}finis sur $J\times \Omega$ et
l'espace $\Lambda(J)$ dont la d\'{e}finition est donn\'{e}e dans
\cite{stoc} p.$24$. On montre suivant \cite{thieu} que l'on peut
calculer les d\'{e}riv\'{e}es $\cal D$ et $\cal D^2$ sur des \'{e}l\'{e}ments de
cet espace (\textit{cf} \cite{stoc} p.$26$). On rappelle le
\begin{theoreme}
Soit $X\in \Lambda$ et $f\in C^{1,2}(I\times \R^d)$ telle que
$\partial_t f$, $\nabla f$ and $\partial_{ij}f$ sont born\'{e}es. On
obtient en adoptant la convention d'Einstein sur la sommation des
indices
\begin{equation}
(\cal{D}X_t)_k =
\left(b-\frac{1}{2p_t}\partial_j(a^{kj}p_t)+\frac{i}{2p_t}\partial_j(a^{kj}p_t)\right)(t,X_t),
\ \mathcal{D} f(t,X_t) = \left(\partial_t f + \mathcal{D} X_t\cdot
\nabla f +\frac{i}{2}a^{kj}\partial_{kj}f\right)(t,X_t)
.\label{deriv_fonc}
\end{equation}
\end{theoreme}

Nous donnons les d\'{e}finitions n\'{e}cessaires \`{a} l'\'{e}tablissement d'un
th\'{e}or\`{e}me de Noether stochastique.

\begin{defi}
Soit $\phi :\R^d \rightarrow \R^d$ un diff\'{e}omorphisme. La
suspension stochastique de $\phi$ est l'application $\Phi :{\cal
P}\rightarrow {\cal P}$ d\'{e}finie par $\forall X\in {\cal P} , \
\Phi(X)_t (\omega )=\phi (X_t (\omega ))$. Dans la suite on notera
indiff\'{e}remment le diff\'{e}omorphisme et sa suspension.\\
De plus un groupe \`{a} un param\`{e}tre de transformations $\Phi_s :
\Upsilon\rightarrow \Upsilon$, $s\in \R$, o\`{u} $\Upsilon \subset
{\cal P}$, est appel\'{e} un groupe $\phi$-suspendu agissant sur
$\Upsilon$ s'il existe un groupe \`{a} un param\`{e}tre de diff\'{e}omorphisme
$\phi_s :\R^d \rightarrow \R^d$, $s\in \R$, tel que pour tout
$s\in \R$, $\Phi_s$ soit une suspension stochastique de $\phi_s$,
et pour tout $X\in \Upsilon$, $\Phi_s (X)\in \Upsilon$.
\end{defi}

\begin{defi}
Un groupe \`{a} un param\`{e}tre de diff\'{e}omorphismes est dit admissible si
$\Phi=\{\phi_s\}_{s\in\R}$ est un groupe \`{a} un param\`{e}tre de
$C^2$-diffeomorphismes sur $\R^d$ tel que $(s,x)\mapsto \partial_
x\phi_s(x)\ \mbox{\rm est}\ C^2$ et tel que la formule
(\ref{deriv_fonc}) reste vrai pour toute fonction $\phi_s$ du
groupe.
\end{defi}

La derni\`{e}re condition peut para\^{\i}tre tr\`{e}s restrictive, mais elle
est v\'{e}rifi\'{e}e pour les groupes \`{a} un param\`{e}tre de diff\'{e}omorphismes
affines de $\R^d$, ce qui est important dans le cas classique
(\cite{ar} p.$89-90$).

\begin{lemm}\label{tecnicnoether}
Soit $\Phi=(\phi_s)_{s\in\R}$ une suspension stochastique d'un
groupe admissible \`{a} un param\`{e}tre de diff\'{e}omorphismes. Alors pour
tout $X\in\Lambda$, et pour tout $(t,s)\in I\times\R$
l'application $s\mapsto \cal{D}(\Phi_s X)_t$ est de classe $C^1$
$p.s.$ et $\d \partial_s[ \mathcal{D}(\phi_s(X))]=\mathcal{D}
\left[
\partial_s \phi_s(X)\right]$ $p.s.$ .
\end{lemm}

\begin{defi}
Soit $\Phi=(\phi_s)_{s\in \R}$ ne suspension stochastique d'un
groupe admissible \`{a} un param\`{e}tre de diff\'{e}omorphismes et $L: \Cu
(I)\rightarrow \Cu_{\C} (I)$. La fonctionnelle $L$ est invariante
sous $\Phi$ si pour tout $s\in\R$ et $X\in\Cu(J)$,  $L(\phi_s
X,\cal{D}(\phi_s(X)))=L(X,\cal{D}X) $.
\end{defi}

Le th\'{e}or\`{e}me s'\'{e}nonce alors :
\begin{theoreme}
\label{noether} Soit $F_J$ la fonctionnelle d\'{e}finie sur
$\Xi\cap\Lambda(J)$ par $\d F_J (X)=E\left [\int_J L(X_t,{\cal D}
X_t ) dt \right ]$, o\`{u} $L$ est un lagrangien admissible invariant
sous le groupe admissible \`{a} un param\`{e}tre de diff\'{e}omorphisme
$\Phi=(\phi_s)_{s\in\R}$. Soit $X^0\in \Xi\cap\Lambda(J)$ un point
$\Cu(J)$-stationnaire de $F_J$. On pose $ Y_t(s) =\Phi_s(X^0)_t$.
Alors $\d \frac{d}{dt}E\left[\partial_v L(X^0,\cal{D}X^0_t) \cdot
\frac{\partial Y_t}{\partial s}(0)\right]=0$.
\end{theoreme}

\textbf{D\'{e}monstration.} On pose $V_t(s)=(Y_t(s),{\cal D} Y_t(s))$.
Comme $L$ est invariant sous $\Phi=\{ \phi_s \}_{s\in \R}$, on a
$\d {\partial \over
\partial s} L(V_t(s))=0 \quad (p.s.).$ Comme pour
tout $t\in J$ et tout $\omega\in\Omega$, $Y_t(\cdot)(\omega)\in
C^1(\R)$ et ${\cal D} Y_t(\cdot)(\omega)\in C^1(\R)$, on obtient
$\d
\partial_x L(V_t(s)) \cdot \d {\partial Y_t \over
\partial s} +\partial_v L(V_t(s)) \cdot {\partial {\cal D} Y_t \over
\partial s} =0 \quad (p.s.)$. En utilisant le lemme (\ref{tecnicnoether}), cette \'{e}quation est \'{e}quivalente \`{a} $\d
\partial_x L(V_t(s)) \cdot \d {\partial Y_t \over \partial s} +\partial_v L(V_t(s))
\cdot {\cal D} \left ( {\partial Y_t \over \partial s}\right ) =0
\quad (p.s.)$. Comme $X^0=Y(0)$ est un $\Cu(J)-$point stationaire
de $F_J$, on a $\d
\partial_x L(V_t(0)) =\overline{{\cal D}} \left[\partial_v L(V_t(0))\right]$.
On en d\'{e}duit alors $\d \overline{{\cal D}} \left[\partial_v
L(V_t(0))\right] \cdot {\partial Y_t \over \partial s}(0)
+\partial_v L(V_t(0)) \cdot {\cal D} \left ( {\partial Y_t \over
\partial s}(0)\right ) =0 \quad (p.s.)$. D'o\`{u} $\d E\left [ \overline{{\cal D}} \left[\partial_v L(V_t(0))\right] \cdot {\partial Y \over \partial s}(0)
+\partial_v L(V_t(0)) \cdot {\cal D} \left ( {\partial Y_t \over
\partial s}(0)\right ) \right] =0$.\\
Avec le lemme (\ref{formule2}), il vient $ \d {d\over dt} E\left [
\partial_v L(V_t(0)) \cdot {\partial Y_t \over
\partial s}(0) \right ] =0$. $\square$

Ce th\'{e}or\`{e}me sugg\`{e}re d'introduire la notion suivante d'int\'{e}grale
premi\`{e}re :
\begin{defi}
Soit $L$ un lagrangien admissible. Une fonctionnelle $I:
L^2(\Omega) \rightarrow \C$ est une int\'{e}grale premi\`{e}re pour
l'\'{e}quation d' Euler-Lagrange stochastis\'{e}e associ\'{e}e \`{a} $L$ si $ \d
{d\over dt} \left[ I(X_t) \right ] =0,$ pour tout $X$ satisfaisant
une \'{e}quation d'Euler-Lagrange stochastique.
\end{defi}

% Les remerciements sont dans une section, sans num\`{E}rotation

%\section*{Remerciements}
% Remerciements - texte ici

\end{document}